\documentclass[11pt]{amsart}
\usepackage[english]{babel}
\usepackage{amscd,amssymb,amsmath}
\usepackage{graphicx}

\def\abd{\mathop {\rm \leftharpoondown \!\!\!\!\! \rightharpoondown} }
\begin{document}
\title[Discrete analogue of the operator]{Construction of discrete analogue of the differential operator
$\frac{d^4}{dx^4}+2\frac{d^2}{dx^2}+1$ and its properties}%
\author{A.R. Hayotov}%
\address{Institute of Mathematics, National University of Uzbekistan, Do`rmon yo`li str. 29,  Tashkent-100125, Uzbekistan}%
\email{hayotov@mail.ru}%

\subjclass{65D32}%
\keywords{Discrete argument function, discrete analogue of the diffential operator, differential operator, Fourier transformation}%

\begin{abstract}
In the present paper the discrete analogue of the differential operator $\frac{d^4}{dx^4}+2\frac{d^2}{dx^2}+1$
is constructed and its some properties are proved.
\end{abstract}
\maketitle

The optimization problem of approximate integration formulas in the modern sense appears as the problem
of finding the minimum of the norm of a error functional $\ell$ given on some set of functions.

The minimization problem of the norm of the error functional by coefficients was reduced in [1]
to the system of difference equations of Wiener-Hopf type in the space $L_2^{(m)}$, where $L_2^{(m)}$
is the space of functions with square integrable $m-$th generalized derivative. Existence and unique\-ness
of the solution of this system was proved by S.L. Sobolev. In the work [1] the description
of some analytic algorithm for finding the coefficients of optimal formulas is given. For this S.L. Sobolev defined and investigated
the discrete analogue of the differential operator
of the poly\-har\-monic  operator $\Delta^m$.
The problem of construction of the discrete opera\-tor $D_{hH}^{(m)}[\beta]$ for  $n-$ dimensional case was very hard. In one dimensional
case the discrete analogue $D_{h}^{(m)}[\beta]$ of the differential operator $\frac{d^{2m}}{dx^{2m}}$ was constructed by
Z.Zh. Zhamalov and Kh.M. Shadimetov [2,3].

Further, in the work [4] the discrete analogue of the differential  $\frac{d^{2m}}{dx^{2m}}-\frac{d^{2m-2}}{dx^{2m-2}}$ was constructed.
The constructed discrete analogue of the operator $\frac{d^{2m}}{dx^{2m}}-\frac{d^{2m-2}}{dx^{2m-2}}$ was applied for finding the coefficients of
optimal quadrature formulas in the space $W_2^{(m,m-1)}(0,1)$
(see [5]).

Here, we  mainly
use a concept of functions of a discrete argument and the corresponding operations
(see [1]).  For completeness we give some of
definitions.

Let $[\beta]=h\beta$, $\beta\in \mathbb{Z}$, $h=\frac{1}{N}$, $N=1,2,...$ .
Assume that $\varphi (x)$ and $\psi (x)$ are real-valued functions
defined on the real line $\mathbb{R}$.
\smallskip

{\sc Definition 1.} The function $\varphi[\beta]$ is a
\emph{function of discrete argument} if it is given on some set
of integer values of $\beta $.
\smallskip

 {\sc Definition 2.} The \emph{inner product} of two
discrete argument functions $\varphi[\beta]$ and $\psi[\beta]$
is given by
\[
\left[ {\varphi ,\psi } \right] = \sum_{\nu  =  - \infty }^\infty
\varphi[\beta] \cdot \psi[\beta],\] if the series on the right
hand side  converges absolutely.
\smallskip

{\sc Definition 3.} The \textit{convolution} of two functions
$\varphi[\beta]$ and $\psi[\beta]$ is  the inner
product
\[\varphi[\beta]*\psi[\beta] = \left[ {\varphi [\gamma],\psi[\beta  -\gamma]} \right] = \sum_{\gamma  =  -
\infty }^\infty  {\varphi [\gamma] \cdot \psi [\beta  - \gamma]}.
\]

In the present paper we consider the problem of construction of the discrete function $D[\beta]$ which satisfies the following equation
$$
D[\beta]*\psi[\beta]=\delta[\beta], \eqno (1)
$$
where
$$
G[\beta]=\frac{\mathrm{sign}[\beta]}{4}(\sin[\beta]-[\beta]\cdot
\cos[\beta]), \eqno (2)
$$
$\delta[\beta]=\left\{
\begin{array}{ll}
1, & \beta=0,\\
0, &\beta\neq 0
\end{array}
\right. $ is the discrete delta function.

The discrete function $D[\beta]$ has important role in calculation of the coef\-fici\-ents of the optimal quadrature
formulas in the Hilbert space
\[
K_2(P_2)=\Bigl\{\varphi:[0,1]\to \mathbb{R}\ \Bigm|\ \varphi' \mbox{ is
absolutely continuous and } \varphi''\in L_2(0,1)\Bigr\},
\]
equipped with the norm
$$
\|\varphi\|=\left\{\int\limits_0^1\left(P_2\left(d/dx\right)
\varphi(x)\right)^2dx\right\}^{\frac12},
$$
where $P_2(d/dx)=d^2/dx^2+1$.

The equation (1) is the discrete analogue of the following equation
$$
\left(\frac{d^4}{dx^4}+2\frac{d^2}{dx^2}+1\right)G(x)=\delta(x),
\eqno (3)
$$
where $G(x)=\frac{\mathrm{sign}(x)}{4}(\sin x-x\cdot \cos x)$,
$\delta(x)$ is Dirac's delta function. Moreover the discrete function
$D[\beta]$ has similar properties as the differential operator
$\frac{d^4}{dx^4}+2\frac{d^2}{dx^2}+1$, i.e. the zeros of the discrete operator $D[\beta]$ are the discrete functions
corresponding to the
zeros of the operator $\frac{d^4}{dx^4}+2\frac{d^2}{dx^2}+1$.
The discrete function $D[\beta]$ is called \emph{the discrete analogue}  of the differential operator
$\frac{d^4}{dx^4}+2\frac{d^2}{dx^2}+1$.

The main results of the present work are the following theorems

\textbf{Theorem 1.} {\it The discrete analogue of the differential operator  $\frac{d^4}{dx^4}+2\frac{d^2}{dx^2}+1$ satisfying
equation (1) has the form
$$
D[\beta]=\frac{2}{\sin h-h\cdot \cos h}\left\{
\begin{array}{ll}
A_1\cdot \lambda_1^{|\beta|-1},& |\beta|\geq 2,\\
1+A_1,& |\beta|=1,\\
\frac{2h\cos 2h-\sin 2h}{\sin h-h\cos h}+\frac{A_1}{\lambda_1},&
\beta=0,
\end{array}
\right.\eqno (4)
$$
where $$ A_1=\frac{4h^2\sin^4h\cdot
\lambda_1^2}{(\lambda_1^2-1)(\sin h-h\cos h)^2},
\ \lambda_1=\frac{2h-\sin 2h-2\sin h\cdot \sqrt{h^2-\sin^2h}}{2(h\cos
h-\sin h)},
$$
$|\lambda_1|<1,\ h$ is a small parameter.}

\textbf{Theorem 2.} {\it The discrete analogue  $D[\beta]$
of the differential operator $\frac{d^4}{dx^4}+2\frac{d^2}{dx^2}+1$ satisfies
the following equalities

1) $D[\beta]*\sin[\beta]=0,$

2) $D[\beta]*\cos[\beta]=0,$

3) $D[\beta]*[\beta]\sin[\beta]=0,$

4) $D[\beta]*[\beta]\cos[\beta]=0,$

5) $D[\beta]*\psi[\beta]=\delta[\beta].$

Here $G[\beta]$ is defined by (2), and $\delta[\beta]$
is discrete delta function.}\\[0.5mm]

In the proofs of these theorems we need the following well known formulas from the theory of
generalized (distribution) functions and Fourier transformations (see, for instance, [1])
$$
F[\varphi(x)]=\int\limits_{-\infty}^{\infty}\varphi(x)e^{2\pi
ipx}dx,\ \
F^{-1}[\varphi(p)]=\int\limits_{-\infty}^{\infty}\varphi(p)e^{-2\pi
ipx}dp,\eqno (5)
$$
$$
F[\varphi*\psi]=F[\varphi]\cdot F[\psi], \eqno (6)
$$
$$
F[\varphi\cdot \psi]=F[\varphi]* F[\psi], \eqno (7)
$$
$$
F[\delta^{(\alpha)}(x)]=(-2\pi ip)^{\alpha},\ \ F[\delta(x)]=1,
\eqno(8)
$$
$$
\delta(hx)=h^{-1}\delta(x),\eqno (9)
$$
$$
\delta(x-a)\cdot f(x)=\delta(x-a)\cdot f(a), \eqno (10)
$$
$$
\delta^{(\alpha)}(x)*f(x)=f^{(\alpha)}(x), \eqno (11)
$$
$$
\phi_0(x)=\sum\limits_{\beta=-\infty}^{\infty}\delta(x-\beta),\ \
\ \sum\limits_{\beta}e^{2\pi i
x\beta}=\sum\limits_{\beta}\delta(x-\beta). \eqno (12)
$$

 \textbf{Proof of Theorem 1.}

According to the theory of periodic generalized functions and Fourier transformations
instead of the function $D[\beta]$ it is convenient to search the harrow-shaped function
(see [1])
$$
\stackrel{\abd}{D}(x)=\sum\limits_{\beta=-\infty}^{\infty}
D[\beta]\delta(x-h\beta).
$$

In the class of harrow-shaped functions equation (1) will be in the following form
$$
\stackrel{\abd}{D}(x)*\stackrel{\abd}{G}(x)=\delta(x), \eqno
(13)
$$
where $\stackrel{\abd}{G}(x)=\sum\limits_{\beta=-\infty}^{\infty}
G[\beta]\delta(x-h\beta)$ is the harrow-shaped function corresponding to the function
$G[\beta]$.

Applying the Fourier transformation to both sides of equation (13)  and taking into account
(6), (8) we have
$$
F[\stackrel{\abd}{D}(x)]=1/F[\stackrel{\abd}{G}(x)].
\eqno (14)
$$
Now we calculate the Fourier transformation Òåïåðü
$F[\stackrel{\abd}{G}(x)]$. Using equalities (10) and (12),
we get
$$
\stackrel{\abd}{G}(x)=h^{-1}G(x)\cdot \phi_0(h^{-1}x). \eqno
(15)
$$
Further, use of equation (9) gives
$$
F[\phi_0(h^{-1}x)]=h\phi_0(hp). \eqno (16)
$$
Then, taking account of (15), (16) and (7), we obtain
$$
F[\stackrel{\abd}{G}(x)]=F[G(x)]*\phi_0(hp). \eqno (17)
$$
To calculate the Fourier transformation $F[G(x)]$
we use equation (3). Taking into account (11), we rewrite equation (3) in the following form
$$
(\delta^{(4)}(x)+2\delta^{(2)}(x)+\delta(x))*G(x)=\delta(x).
$$
Hence, keeping in mind (6), (8), we have
$$
F[G(x)]=\frac{1}{(2\pi p-1)^2(2\pi p+1)^2}. \eqno (18)
$$
Taking into account (18) from (17) we get
$$
F[\stackrel{\abd}{G}(x)]=\frac{h^3}{(2\pi)^4}
\sum\limits_{\beta=-\infty}^{\infty}\frac{1}{(\beta-h(p+\frac{1}{2\pi}))^2
(\beta-h(p-\frac{1}{2\pi}))^2}.
$$
Then from (14) we have
$$
F[\stackrel{\abd}{D}](p)=\frac{(2\pi)^4}{h^3}
\left[\sum\limits_{\beta=-\infty}^{\infty}\frac{1}{(\beta-h(p+\frac{1}{2\pi}))^2
(\beta-h(p-\frac{1}{2\pi}))^2}\right]^{-1}. \eqno (19)
$$
Suppose the Fourier series of the function  $F[\stackrel{\abd}{D}](p)$ has the following form
$$
F[\stackrel{\abd}{D}](p)=\sum\limits_{\beta=-\infty}^{\infty}\hat
D[\beta]e^{2\pi i ph\beta},\eqno (20)
$$
where $\hat D[\beta]$ is the Fourier coefficients of the function
$F[\stackrel{\abd}{D}](p)$, i.e.
$$
\hat D[\beta]=\int\limits_0^{h^{-1}}F[\stackrel{\abd}{D}](p)\
e^{-2\pi i ph\beta}dp. \eqno (21)
$$
Applying the inverse Fourier transformation to both sides of (20)  we obtain the following
harrow-shaped function
$$
\stackrel{\abd}{D}(x)=\sum\limits_{\beta=-\infty}^{\infty}\hat
D[\beta]\delta(x-h\beta).
$$
Then according to the definition of harrow-shaped functions the discrete function
$\hat D[\beta]$ is searching function of discrete argument
$D[\beta]$. Here for finding of the function   $\hat D[\beta]$ we will not use the formula (21).
We will find it by the following way.

To calculate the series (19) we use the following well known formula from
the residual theory (see [6], p.296)
$$
\sum\limits_{\beta=-\infty}^{\infty}f(\beta)=-\sum\limits_{z_1,z_2,...,z_n}\mbox{
res}(\pi \mathrm{ctg}(\pi z)\cdot f(z)),\eqno (22)
$$
where $z_1,z_2,...,z_n$ are poles of the function $f(z)$.

We denote
$f(z)=\frac{1}{(z-h(p+\frac{1}{2\pi}))^2(z-h(p-\frac{1}{2\pi}))^2}$.
Here $z_1=h(p+\frac{1}{2\pi})$ and $z_2=h(p-\frac{1}{2\pi})$ are the poles of order 2. Then taking into account (22) from (19) we have
$$
F[\stackrel{\abd}{D}](p)=-\frac{(2\pi)^4}{h^3} \left[
\sum\limits_{z_1,z_2}\mbox{ res}(\pi \mathrm{ctg}(\pi z)\cdot
f(z)) \right]^{-1}. \eqno (23)
$$
Since
$$
\mathop{\mbox{res}}\limits_{z=z_1}(\pi \mathrm{ctg}(\pi z)\cdot
f(z))=-\frac{2\pi^4}{h^3}\left(\frac{h}{1-\cos(2\pi
hp+h)}+\frac{\sin(2\pi hp+h)}{1-\cos(2\pi hp+h)}\right),
$$
$$
\mathop{\mbox{res}}\limits_{z=z_2}(\pi \mathrm{ctg}(\pi z)\cdot
f(z))=-\frac{2\pi^4}{h^3}\left(\frac{h}{1-\cos(2\pi
hp-h)}-\frac{\sin(2\pi hp-h)}{1-\cos(2\pi hp-h)}\right).
$$
Denoting $\lambda=e^{2\pi iph}$, using the last two equalities and taking into account the following well known
formulas
$$
\cos z=\frac{e^{zi}+e^{-zi}}{2},\ \ \sin
z=\frac{e^{zi}-e^{-zi}}{2i}
$$
after some simplifications from (23) for $F[\stackrel{\abd}{D}](p)$
we get
$$
F[\stackrel{\abd}{D}](p)= \frac{2}{\sin h-h\cos h}\cdot \frac{
\lambda^4-4\cos h\ \lambda^3+(2\cos(2h)+4)\lambda^2-4\cos h\
\lambda+1}{\lambda\left(\lambda^2+
\frac{2h-\sin(2h)}{sinh-h\cosh}\lambda+1\right)}. \eqno (24)
$$

To find the explicit form of the discrete operator  $D[\beta]$
the equality (24) we expand to the sum of elementary fractions. Since the polynomial
$Q_2(\lambda)=\lambda^2+ \frac{2h-\sin(2h)}{sinh-h\cosh}\lambda+1$
has two real roots
$$
\lambda_1=\frac{2h-\sin(2h)-2\sin h\cdot \sqrt{h^2-\sin^2h}}
{2(h\cos h-\sin h)},
$$
$$
\lambda_2=\frac{2h-\sin(2h)+2\sin h\cdot \sqrt{h^2-\sin^2h}}
{2(h\cos h-\sin h)}
$$
and $\lambda_1\cdot\lambda_2=1$, $|\lambda_1|<1$.\\
Then from  (24) we have
$$
\frac{2}{\sin h-h\cos h}\cdot \frac{ \lambda^4-4\cos h\
\lambda^3+(2\cos(2h)+4)\lambda^2-4\cos h\
\lambda+1}{\lambda\left(\lambda^2+ \frac{2h-\sin(2h)}{\sin
h-h\cosh}\lambda+1\right)}= $$ $$= \frac{2}{\sin h-h\cos h}\cdot
\left(\lambda+\frac{2h\cos(2h)-\sin(2h)}{\sin h-h\cos
h}+\frac{A}{\lambda}+\frac{A_1}{\lambda-\lambda_1}+\frac{B_1}{\lambda-\lambda_2}\right).
\eqno (25)
$$
For finding unknown coefficients $A$, $A_1$ and  $B_1$
in the equation (25) we put $\lambda=0$,
$\lambda=\lambda_1$ and $\lambda=\lambda_2$.\\
For $\lambda=0$
$$
A=1,\eqno (26)
$$
for $\lambda=\lambda_1$
$$
A_1=\frac{\lambda_1^4-4\cos h\
\lambda_1^3+(2\cos(2h)+4)\lambda_1^2-4\cos h\
\lambda_1+1}{\lambda_1^2-1},\eqno (27)
$$
for $\lambda=\lambda_2$
$$
B_1=\frac{\lambda_2^4-4\cos h\
\lambda_2^3+(2\cos(2h)+4)\lambda_2^2-4\cos h\
\lambda_2+1}{\lambda_2^2-1}.
$$
Hence, taking into account $\lambda_1\cdot \lambda_2=1$, we have
$$
B_1=-\frac{1}{\lambda_1^2}\cdot A_1.
\eqno (28)
$$
Finally, taking account of (26), (27), (28) and  $|\lambda|=1$,
$|\lambda_1|<1$, using the formula for geometric progression from (25) consequently we get {\small
$$
F[\stackrel{\abd}{D}](p)=
$$
$$ = \frac{2}{\sin h-h\cos h}\cdot
\left(\lambda+\frac{2h\cos(2h)-\sin(2h)}{\sin h-h\cos
h}+\frac{A}{\lambda}+\frac{A_1}{\lambda}\frac{1}{1-\frac{\lambda_1}{\lambda}}
-\frac{B_1}{\lambda_2}\frac{1}{1-\frac{\lambda}{\lambda_2}}\right)=
$$
$$ = \frac{2}{\sin h-h\cos h}\cdot
\Bigg(\lambda+\frac{2h\cos(2h)-\sin(2h)}{\sin h-h\cos
h}+\frac{A}{\lambda}+
$$
$$
+ \frac{A_1}{\lambda}\sum\limits_{\gamma=0}^{\infty}
\left(\frac{\lambda_1}{\lambda}\right)^\gamma-
{B_1\lambda_1}\sum\limits_{\gamma=0}^{\infty}
\left(\lambda_1\lambda\right)^\gamma\Bigg)=
$$
$$ = \frac{2}{\sin h-h\cos h}\cdot
\Bigg(\lambda(1+A_1)+\frac{2h\cos(2h)-\sin(2h)}{\sin h-h\cos
h}+\frac{A_1}{\lambda_1}+
$$
$$
+ (1+A_1)\frac{1}{\lambda}+A_1\sum\limits_{\gamma=-2}^{\infty}
\lambda_1^{-\gamma-1}\lambda^\gamma+A_1
\sum\limits_{\gamma=2}^{\infty} \lambda_1^{\gamma-1}\lambda^\gamma
\Bigg)=\sum\limits_{\gamma=-\infty}^{\infty}D[\gamma]\lambda^\gamma.
$$
 }
Hence keeping in mind  $\lambda=e^{2\pi iph}$ we obtain the explicit form
(4) of the discrete function $D[\beta]$.

Theorem 1 is proved.\\[1mm]

Theorem 2 is proved using Definition 3  and by direct calculation
of the left hand sides of 1)-5).\\[1mm]

\textbf{Remark.} \emph{From (4) we note that  $D[\beta]$ is the even function, i.e.
$D[\beta]=D[-\beta]$ and since $|\lambda_1|<1$, then the function
$D[\beta]$ is exponentially decreased as $\beta\to \infty$.}

\end{document}